\documentclass[reqno, 12pt]{amsart}

\usepackage{amssymb}

\usepackage{comment}
\usepackage[left=2.2cm,top=2.4cm,right=2.2cm]{geometry}

\usepackage[hidelinks]{hyperref}

\usepackage[latin1]{inputenc}
\usepackage[initials, alphabetic,backrefs,msc-links]{amsrefs}
\usepackage{mathrsfs}
\usepackage{soul}
\usepackage{amssymb}
\usepackage{amsmath}
\usepackage{amscd,amsthm,amsxtra, array,ulem}

\usepackage{fancybox,color}
\usepackage{enumerate}

\usepackage[dvipsnames]{xcolor}

\newcommand{\N}{\mathbb{N}}

\newcommand{\R}{\mathbb{R}}

\newcommand{\norm}[1]{\left\|#1\right\|}

\numberwithin{equation}{section}

\allowdisplaybreaks

\def\vint_#1{\mathchoice%
        {\mathop{\kern 0.2em\vrule width 0.6em height 0.69678ex depth -0.58065ex
                \kern -0.8em \intop}\nolimits_{\kern -0.4em#1}}%
        {\mathop{\kern 0.1em\vrule width 0.5em height 0.69678ex depth -0.60387ex
                \kern -0.6em \intop}\nolimits_{#1}}%
        {\mathop{\kern 0.1em\vrule width 0.5em height 0.69678ex
            depth -0.60387ex
                \kern -0.6em \intop}\nolimits_{#1}}%
        {\mathop{\kern 0.1em\vrule width 0.5em height 0.69678ex depth -0.60387ex
                \kern -0.6em \intop}\nolimits_{#1}}}
\def\vintslides_#1{\mathchoice%
        {\mathop{\kern 0.1em\vrule width 0.5em height 0.697ex depth -0.581ex
                \kern -0.6em \intop}\nolimits_{\kern -0.4em#1}}%
        {\mathop{\kern 0.1em\vrule width 0.3em height 0.697ex depth -0.604ex
                \kern -0.4em \intop}\nolimits_{#1}}%
        {\mathop{\kern 0.1em\vrule width 0.3em height 0.697ex depth -0.604ex
                \kern -0.4em \intop}\nolimits_{#1}}%
        {\mathop{\kern 0.1em\vrule width 0.3em height 0.697ex depth -0.604ex
                \kern -0.4em \intop}\nolimits_{#1}}}

\newcommand{\aveint}[2]{\mathchoice%
        {\mathop{\kern 0.2em\vrule width 0.6em height 0.69678ex depth -0.58065ex
                \kern -0.8em \intop}\nolimits_{\kern -0.45em#1}^{#2}}%
        {\mathop{\kern 0.1em\vrule width 0.5em height 0.69678ex depth -0.60387ex
                \kern -0.6em \intop}\nolimits_{#1}^{#2}}%
        {\mathop{\kern 0.1em\vrule width 0.5em height 0.69678ex depth -0.60387ex
                \kern -0.6em \intop}\nolimits_{#1}^{#2}}%
        {\mathop{\kern 0.1em\vrule width 0.5em height 0.69678ex depth -0.60387ex
                \kern -0.6em \intop}\nolimits_{#1}^{#2}}}

\def\XXint#1#2#3{{\setbox0=\hbox{$#1{#2#3}{\int}$}
    \vcenter{\hbox{$#2#3$}}\kern-.5\wd0}}

\newcommand{\vertiii}[1]{{\left\vert\kern-0.25ex\left\vert\kern-0.25ex\left\vert #1 
    \right\vert\kern-0.25ex\right\vert\kern-0.25ex\right\vert}}

\newcommand{\vertii}[1]{{\left\vert\kern-0.25ex\left\vert\kern-0.25ex  #1 
    \kern-0.25ex\right\vert\kern-0.25ex\right\vert}}

\makeatletter
\@namedef{subjclassname@2020}{\textup{2020} Mathematics Subject Classification}
\makeatother

\newtheorem{theorem}{Theorem}[section]
\newtheorem{lemma}[theorem]{Lemma}
\newtheorem{corollary}[theorem]{Corollary}
\newtheorem{proposition}[theorem]{Proposition}

\theoremstyle{definition}
\newtheorem{definition}[theorem]{Definition}

\theoremstyle{remark}
\newtheorem{remark}[theorem]{Remark}

\numberwithin{equation}{section}

\begin{document}

\title{
Degenerate Poincaré-Sobolev inequalities via fractional integration
}
\author[A. Claros]{Alejandro Claros}

\address{BCAM -- Basque Center for Applied Mathematics, Bilbao, Spain}
\email{aclaros@bcamath.org}

\address{Universidad del País Vasco / Euskal Herriko Unibertsitatea (UPV/EHU), Bilbao, Spain}
\email{aclaros003@ikasle.ehu.eus}

\thanks{This research is supported by the Basque Government through the BERC 2022-2025 program, by the Ministry of Science and Innovation through Grant PRE2021-099091 funded by BCAM Severo Ochoa accreditation CEX2021-001142-S/MICIN/AEI/10.13039/501100011033 and by ESF+, and by the project PID2023-146646NB-I00 funded by MICIU/AEI/10.13039/501100011033 and by ESF+.}





\begin{abstract}
	We present a local weighted estimate for the Riesz potential in $\R^n$, which improves the main theorem of Alberico, Cianchi, and Sbordone [C. R. Math. Acad. Sci. Paris \textbf{347} (2009)] in several ways. As a consequence, we derive weighted Poincaré-Sobolev inequalities with sharp dependence on the constants. We answer positively to a conjecture proposed by Pérez and Rela [Trans. Amer. Math. Soc. \textbf{372} (2019)] related to the sharp exponent in the $A_1$  constant in the $(p^*,p)$ Poincaré-Sobolev inequality with $A_1$ weights. Our approach is versatile enough to prove Poincaré-Sobolev inequalities for high-order derivatives and fractional Poincaré-Sobolev inequalities with the BBM extra gain factor $(1-\delta)^\frac{1}{p}$.   In particular, we improve one of the main results from Hurri-Syrjänen, Martínez-Perales,  Pérez, and Vähäkangas [Int. Math. Res. Not. \textbf{20} (2023)].
\end{abstract}

\maketitle

\section{Introduction}

Weighted Poincaré-Sobolev inequalities play a key role in proving the local Hölder regularity of weak solutions of degenerate elliptic partial differential equations. In the celebrated paper \cite{FKS}, E. Fabes, C. Kenig, and R. Serapioni considered the following degenerate elliptic PDE, namely the operator $Lu=\operatorname{div}(A(x)\nabla u)$, where $A$ is an $n\times n$ symmetric matrix in $\R^n$ satisfying the following degenerate ellipticity condition: there exist $C_1, C_2>0$ such that
\begin{equation*}
	C_1 w(x) |\xi |^2 \le \left \langle A(x) \xi , \xi 	\right \rangle \le C_2 w(x) |\xi |^2
\end{equation*}
where $w$ is a weight in the $A_2$ Muckenhoupt class. The Moser iteration method is a classical and powerful technique to prove the local Hölder regularity of the weak solutions of a given elliptic partial differential equation (see \cite{HKM}). The method has two important steps: the $(2,2)$ Poincaré inequality and the  $(2^*, 2)$ Poincaré-Sobolev inequality, where $2^*=\frac{2n}{n-2}>2$ is the classical Sobolev exponent. In order to apply the Moser iteration method to the degenerate setting, in \cite{FKS} weighted Poincaré and Poincaré-Sobolev inequalities are proved, namely
\begin{equation}\label{intro PS}
	\left(\frac{1}{w(B)}\int_B |f(x)-f_B|^q w(x)dx \right)^\frac{1}{q}\le C_w r(B)\left( \frac{1}{w(B)}\int_B |\nabla f(x)|^p w(x)dx\right)^\frac{1}{p},
\end{equation}
for all balls $B$ with radius $r(B)$, where $f_B=\frac{1}{|B|}\int_B f$ is the average of $f$ over the ball $B$, $1\le p \le q <\infty$ and $w$ is a weight in some $A_p$ class of Muckenhoupt (see Section \ref{Section 2} for the exact definition). This type of inequality has attracted the attention of many authors \cite{ CW1, CW2, CF, FPW98,  Hu, CPER, CPER2, SW}.

A general method for obtaining inequalities of this type was introduced in \cite{FPW98}, and then refined in \cite{MP98}. The main result has recently been improved in \cite{CP}. This abstract approach consists of the so-called self-improving property, which satisfies Poincaré inequalities and generalized Poincaré inequalities. The self-improving theory unifies the Poincaré inequalities and the John-Nirenberg theorem from the $\operatorname{BMO}$ space. 

In a recent work, C. Pérez and E. Rela have introduced a new theory of self-improvement of generalized Poincaré inequalities. As a consequence, the weighted Poincaré-Sobolev inequalities obtained in \cite{FKS} are improved. They also pay attention to the question of studying quantitative estimates of the constant $C_w$ in \eqref{intro PS} and obtaining a precise control of $q$ in terms of $n$, $p$, and $w$. Concretely, in Corollary 1.22 of \cite{CPER}, the inequality \eqref{intro PS} is proved for all $w\in A_1$ with $q=\frac{np}{n-p}$ the classical Sobolev exponent and constant $C_w=c_{n, p} [w]_{A_1}$. The authors prove that the exponent for the $A_1$ constant cannot be better than $\frac{1}{p}$ and conjecture that the same inequality holds with the constant $C_w$ replaced by $C_w'=c_{n,p}[w]_{A_1}^{1/p}$ (cf. Conjecture 1.23 of \cite{CPER}). As a consequence of our main result, we answer positively to this conjecture (see Theorem \ref{thm ps 1} and the consequent remark). 

The classical approach for obtaining Poincaré-Sobolev inequalities is through the use of a well-known subrepresentation formula of $f$ in terms of a fractional integral of the gradient; namely, there exists $c_n>0$ such that for each ball $B$,
\begin{equation}\label{subrepresentation formula 1}
	|f(x)-f_B|\le c_n I_1(|\nabla f|\chi_B)(x), \qquad x\in B,
\end{equation}
where $I_1$ is the fractional integral operator of order $1$, also called Riesz potential. The Riesz potential of a function $f: \R^n\longrightarrow \R$ is  defined by 
\begin{equation*}
	I_\alpha f(x)= \int_{\R^n} \frac{f(y)}{|x-y|^{n-\alpha}}dy,
\end{equation*}
for each $x\in \R^n$, where $\alpha\in (0,n)$. The classical theory states that $I_\alpha : L^p \longrightarrow L^q$ if $1< p <\frac{n}{\alpha}$, where $q$ is defined by $\frac{1}{p}-\frac{1}{q}=\frac{\alpha}{n}$, and in the case $p=1$ we have $I_\alpha : L^1 \longrightarrow L^{\frac{n}{n-\alpha}, \infty}$ (see for instance the book \cite{KLV1}). In \cite{MW}, Muckenhoupt and Wheeden studied weighted bounds for the fractional integral operator (see also \cite{LMPT} and the references therein for a modern approach).

In the present paper, we are interested in studying local bounds of $I_\alpha$ with $A_p$ weights. This framework for obtaining weighted Poincaré-Sobolev inequalities was first used in \cite{CF}. This type of inequality was studied in \cite{ACS}. We state the main result of that paper.

\begin{theorem}[\cite{ACS}]\label{thm acs}
	Let $n\ge 2$, $\alpha \in (0,n)$ and $1<p<\frac{n}{\alpha}$. Let $w\in A_p$ and given $r\ge 1$ consider $q$ defined by the relation
	\begin{equation*}
		\frac{1}{p}- \frac{1}{q} = \frac{\alpha}{n} \frac{1}{r}.
	\end{equation*}
	Then, there exist positive constants $k=k(p,n)$ and $C=C(\alpha, p, n)$ such that if 
	\begin{equation} \label{q}
		p-\frac{k}{[w]_{A_p}^\frac{1}{p-1}} < r \le p,
	\end{equation}
	then 
	\begin{equation}\label{Eq1}
		\left( \frac{1}{w(B)} \int_{B} |I_\alpha f(x)|^q w(x)dx\right)^\frac{1}{q} \le C [w]_{A_p}^{\frac{nr-\alpha }{nr(p-1)}} r(B)^\alpha  \left( \frac{1}{w(B)} \int_{B} | f(x)|^p w(x)dx\right)^\frac{1}{p},
	\end{equation}
	for any ball $B\subset \R^n$ of radius $r(B)$ and every function $f\in L^p(B, w)$ (continued by $0$ outside $B$). Moreover, the exponent $ \frac{nr-\alpha }{nr(p-1)}$ of the $A_p$  constant is sharp, in the sense that the statement is false if $\frac{nr-\alpha}{nr(p-1)}$ is replaced by any smaller exponent. 
\end{theorem}

Our main result, Theorem \ref{coro1}, improves Theorem \ref{thm acs} in several ways. First, we improve the $r$'s obtained in \eqref{q}, which produces a larger $q$. Second, we obtain a better quantitative bound in \eqref{Eq1}, in terms of a mixed $A_p - A_\infty$ constant (see \cite{HPR}). We also extend the result to the case $p=1$, replacing the $L^q$ average on the left-hand side with a weak-$L^q$ average.  

The method is very versatile. We can derive other types of inequalities using different subrepresentation formulas and our main results. For instance, we can obtain Poincaré-Sobolev inequalities for high-order derivatives.
 In addition, we obtain weighted fractional Poincaré-Sobolev inequalities, namely,
\begin{equation*}
	\left(\frac{1}{w(B)}\int_B |f(x)-f_B|^q w(x)dx \right)^\frac{1}{q}\le C_w r(B) ^\delta \left( \frac{1}{w(B)} \int_B \int_B \frac{|f(x)-f(y)|^p}{|x-y|^{n+\delta p }} dy \,  w(x)dx \right)^\frac{1}{p}. 
\end{equation*}
where the object on the right-hand side is the local weighted Gagliardo seminorm of $f$, $1\le p\le q<\infty$ and $0<\delta<1$. We can prove this inequality with the interesting BBM extra gain factor $(1-\delta)^\frac{1}{p}$ in the inequality (see Theorem \ref{Fractional} and the preceding comments for more details). 

We conclude by answering the problem of finding the optimal exponent $q$  in \eqref{intro PS} for the class $A_p$. Specifically, for a fixed $1\le p<n$ and a fixed $w\in A_p$, we define an exponent $p^*_w$  and prove that the inequality (1.1) cannot hold for exponents $q>p^*_w$ (see Definition \ref{def p*w} and the following results).

Let us summarize the paper by briefly describing the content of each of the sections. In the next section, we will state the results and make further comments and remarks. In Section \ref{Section Ia}, we state the main result about weighted local bounds for $I_\alpha$.  In Section \ref{Section PS}, we obtain weighted Poincaré-Sobolev inequalities and fractional Poincaré-Sobolev inequalities.  In Section \ref{Section PS m}, we derive Poincaré-Sobolev inequalities for high-order derivatives. In Section \ref{Section exp}, we state a result about the best possible exponent $q$ in \eqref{intro PS} for each weight $w\in A_p$. In Section \ref{Section 2}, we collect some definitions and auxiliary results. In Section \ref{Section 3}, we provide the proof of the results stated in Section \ref{Section Ia}, and in Section \ref{Section 4}, we prove the remaining results.

\section{Statement of the results}\label{Section 1.1}

\subsection{Weighted local bounds for fractional integrals with $A_p$ weights}\label{Section Ia}

We begin this section by noting that the $A_p$ condition is necessary for \eqref{Eq1}. Let us assume that for a weight $w$ we have the inequality
	\begin{equation}\label{eq4}
		\left( \frac{1}{w(B)} \int_{B} |I_\alpha f(x)|^{q} w(x)dx\right)^\frac{1}{q} \le C  r(B)^\alpha  \left( \frac{1}{w(B)} \int_{B} | f(x)|^p w(x)dx\right)^\frac{1}{p},
	\end{equation}
	for fixed $p, q$, for any ball $B$ and every $f\in L^p(B,w)$, then necessarily $w\in A_p$.  Observe that for any ball $B$ and any non-negative function $f$, we have
	\begin{equation*}
		\frac{r(B)^\alpha }{|B|} \int_B f(y)dy \le M_\alpha f(x)\le c_n I_\alpha f(x),
	\end{equation*}
	for all $x\in B$, where $M_\alpha$ is the fractional maximal function defined by 
	\begin{equation*}
		M_\alpha f(x)= \sup_B \left(\frac{r(B)^\alpha }{|B|} \int_B |f(y)|dy\right)\chi_B(x).
	\end{equation*}
	Then, applying \eqref{eq4} we have,
	\begin{equation*}
		\frac{r(B)^\alpha }{|B|} \int_B f(y)dy \le C  r(B)^\alpha  \left( \frac{1}{w(B)} \int_{B} | f(x)|^p w(x)dx\right)^\frac{1}{p},
	\end{equation*} 
	and this is equivalent to $w\in A_p$ (see for instance \cite{GCRdF}). We have proved that if \eqref{eq4} holds for a weight $w$, then $w\in A_p$. This observation shows that it is necessary to consider weights in class $A_p$ for this type of inequality.

Our first result is the following theorem, which is an improvement of the main result of \cite{ACS},  Theorem \ref{thm acs} above. 
 
\begin{theorem}\label{coro1}
	Let $n \ge 1$, $\alpha \in (0,n)$ and $1<p<\frac{n}{\alpha}$. Let $w\in A_r$ with $1\le r\le p$. Consider $q_r$ defined by the relation
	\begin{equation}\label{exp1}
		\frac{1}{p}- \frac{1}{q_r} = \frac{\alpha}{n}\frac{\tau_n [\sigma]_{A_\infty} }{1+r(\tau_n [\sigma]_{A_\infty} -1)}, 
	\end{equation}
	where $\sigma= w^{1-r'}$ and $\tau_n=2^{n+1}$ if $r>1$ and $q_1=p^*_\alpha = \frac{n p }{n- \alpha p }$. Then, there exists a dimensional constant $c_n>0$ such that, if $r>1$ we have 
	\begin{equation*}
		\left( \frac{1}{w(B)} \int_{B} |I_\alpha f(x)|^{q_r} w(x)dx\right)^\frac{1}{q_r} \le \frac{c_n  }{\alpha } p^*_\alpha  \left(p'\right)^\frac{1}{q_r}  [w]_{A_r}^{\frac{1}{p}} [\sigma]_{A_\infty}^\frac{1}{q_r}  r(B)^\alpha  \left( \frac{1}{w(B)} \int_{B} | f(x)|^p w(x)dx\right)^\frac{1}{p},
	\end{equation*}
	for any ball $B$ and every function $f\in L^p(B, w)$, and if $r=1$ we have, 
	\begin{equation*}
		\left( \frac{1}{w(B)} \int_{B} |I_\alpha f(x)|^{q_1} w(x)dx\right)^\frac{1}{q_1} \le \frac{c_n  }{\alpha } p^*_\alpha  \left(p'\right)^\frac{1}{q_1}  [w]_{A_1}^{\frac{1}{p}}  r(B)^\alpha  \left( \frac{1}{w(B)} \int_{B} | f(x)|^p w(x)dx\right)^\frac{1}{p},
	\end{equation*}
	for any ball $B$ and every function $f\in L^p(B, w)$.
\end{theorem}

	 The proof relies on a weighted version of Hedberg's inequality (see Lemma \ref{lemma Hedberg} below). We establish a local pointwise estimate for the fractional integral operator $I_{\alpha}$ in terms of the Hardy-Littlewood maximal function $M$. Then, to obtain a strong norm inequality, we use the following weighted bound for the maximal operator 
	\begin{equation*}
		\| Mf\|_{L^{p}(w)} \le c_n \left(\frac{p}{p-r}\right)^\frac{1}{p} [w]_{A_r}^\frac{1}{p} \|f\|_{L^p(w)},
	\end{equation*}
	where $w\in A_r$ with $1\le r<p$ (see \eqref{Maximal Lp} below for more details).  This approach, together with the left-openness property of the class $A_r$, naturally produces a mixed-type constant, which includes the $[\sigma]_{A_{\infty}}^{1/q_r}$ term in the estimate.

	  It should be noted that if we used the good-lambda method from \cite{MW}, the dependency on the weight constant would be worse. That is, the inequality
	\begin{equation*}
		\norm{I_\alpha f}_{L^p(w)}\le C_w \norm{M_\alpha f}_{L^p(w)}
	\end{equation*}
	holds for all $0<p<\infty$ and $w\in A_\infty$, where $C_w= C([w]_{A_\infty})$. We note that in \cite{HJ}, it seems to be implicit that the above inequality holds with a linear dependence on the constant, that is, $C_w=c_{n,p} [w]_{A_\infty}$.  The $M_\alpha$ operator is a more manageable operator, which is why the approach in \cite{FKS} follows this method. Nevertheless, it introduces the constant $[w]_{A_\infty}$ into the estimates, and we want to be precise with the dependence on constants. 
	
	\begin{remark}
		 Note that in the case $r=p$, the exponent $q_p$ defined by \eqref{exp1} is larger than the range of $q$'s obtained in Theorem \ref{thm acs}. We remark the improvement of the $A_p$ constant, using basic properties of $A_p$ weights (see Section \ref{Section 2}) we have,
		\begin{align*}
		[w]_{A_p}^{\frac{1}{p}} [\sigma]_{A_\infty} ^\frac{1}{q} \le & c_n [w]_{A_p}^{\frac{1}{p}} [\sigma]_{A_{p'}} ^\frac{1}{q} =  c_n [w]_{A_p}^{\frac{1}{p}} [w]_{A_{p}} ^\frac{p'-1}{q} = c_n  [w]_{A_p}^{\frac{nr-\alpha }{nr(p-1)}}.
		\end{align*} 
		
		 This observation shows that the exponent $\frac{1}{p}$ in $[w]_{A_r}$ and the exponent $\frac{1}{q_r}$ in $[\sigma]_{A_\infty}$ cannot be improved in general. Specifically, if $r = p$, any improvement of these exponents would imply an improved exponent of the $A_p$ constant in Theorem \ref{thm acs}, contradicting its sharpness \cite[p. 1269]{ACS}.
	\end{remark}

The following result is an analogous estimate to the previous theorem, but considering the weak norm on the left-hand side, this allows us to reach the case $p=1$.

\begin{theorem}\label{Teorema I alpha debil 1}
	Let $n\ge 1$, $\alpha \in (0,n)$ and $1\le p<\frac{n}{\alpha}$. Let $w\in A_r$ with $1\le r \le p$. Consider $q_r$ defined in \eqref{exp1} if $r>1$ and $q_1=p^*_\alpha = \frac{n p }{n- \alpha p }$. Then, there exists a dimensional constant $c_n>0$ such that  
	\begin{equation*}
		\left\|I_\alpha f\right\|_{ L^{q_r ,\infty} \left( B , \frac{w(x)dx}{w(B)} \right) } \le \frac{c_n    }{\alpha } p^*_\alpha [w]_{A_r}^{\frac{1}{p}}  r(B)^\alpha  \left( \frac{1}{w(B)} \int_{B} | f(x)|^p w(x)dx\right)^\frac{1}{p},
	\end{equation*}
	for any ball $B$ and every function $f\in L^p(B, w)$. 
\end{theorem}

This result is interesting in itself because we obtain a better constant than the obtained with the strong norm, and in order to derive Poincaré-Sobolev inequalities, it is enough to prove weak weighted estimates.

\subsection{Applications to Poincaré-Sobolev inequalities}\label{Section PS}

We can prove weighted Poincaré-Sobolev inequalities using the subrepresentation formula \eqref{subrepresentation formula 1}, the previous result and the truncation argument (see \cite{H, KO03}), available in the context of Poincaré-Sobolev inequalities.

\begin{theorem}\label{thm ps 1}
	Let $n\ge 2$ and let $1\le p< n$. Let $w\in A_r$ for some $1\le r \le p$  and consider $q_r$ defined by the relation
	\begin{equation}\label{qr PS eq}
		\frac{1}{p}- \frac{1}{q_r} = \frac{1}{n} \frac{\tau_n [\sigma]_{A_\infty} }{1+r( \tau_n [ \sigma ] _{A_\infty} -1)}
	\end{equation}
	where $\sigma = w^{1-r'}$ and $\tau_n=2^{n+1}$ if $r>1$ and $q_1= p^*= \frac{n p }{n- p }$. Then, there exists a dimensional constant $c_n>0$ such that,
		\begin{equation}\label{eq6}
			\left( \frac{1}{w(B)} \int_{B} |f(x)-f_B|^{q_r} w(x)dx\right)^\frac{1}{q_r}  \le c_n p^* [w]_{A_r}^{\frac{1}{p}} r(B)  \left( \frac{1}{w(B)} \int_{B} | \nabla f(x)|^p w(x)dx\right)^\frac{1}{p}
		\end{equation}
		for any ball $B$, and the exponent $ \frac{1 }{p}$ on the $A_r$ constant is sharp. Moreover, using the truncation method, we can self-improve the left-hand side to the  Lorentz norm $L^{q_r, p}$ (see Section \ref{Section 2} for precise definition),	
		\begin{equation*}
			\left\|f-f_B\right\|_{ L^{q_r, p} \left( B , \frac{w(x)dx}{w(B)} \right) }\le  c_n p^* [w]_{A_r}^{\frac{1}{p}} r(B)  \left( \frac{1}{w(B)} \int_{B} | \nabla f(x)|^p w(x)dx\right)^\frac{1}{p}
		\end{equation*}
	for any ball $B$.
	\end{theorem}

The sharpness of the exponent $ \frac{1 }{p}$ on the $A_r$ constant in \eqref{eq6} is proved in Proposition \ref{sharp1}.

\begin{remark}
	In the case $r=1$, this theorem answers positively to the Conjecture 1.23 from \cite{CPER}. Observe that $q_1=p^*=\frac{n p }{n- p }$ is the classical Sobolev exponent.  In that case, we have proved that there exists a dimensional constant $c_n>0$ such that
	\begin{equation*}
			\left( \frac{1}{w(B)} \int_{B} |f(x)-f_B|^{p^*} w(x)dx\right)^\frac{1}{p^*}  \le c_n p^* [w]_{A_1}^{\frac{1}{p}} r(B)  \left( \frac{1}{w(B)} \int_{B} | \nabla f(x)|^p w(x)dx\right)^\frac{1}{p}
		\end{equation*}
		for all $w\in A_1$, and the exponent $\frac{1}{p}$ in the $A_1$ constant is sharp (cf. Proposition 7.4 \cite{CPER}). 
\end{remark}

\begin{remark}
	  The new weighted Sobolev exponent $q_r$ defined by
	  \begin{equation*}
		\frac{1}{p}- \frac{1}{q_r} = \frac{1}{n} \frac{\tau_n [\sigma]_{A_\infty} }{1+r( \tau_n [ \sigma ] _{A_\infty} -1)}
	\end{equation*} 
	is, to the best of our knowledge, new to the literature. This exponent is larger than the best weighted Sobolev exponents obtained in \cite{FPW98} and  \cite{CPER} for the $A_r$ class of weights with $1< r\le p$.
\end{remark}

Now, we focus on Poincaré-Sobolev inequalities with the local Gagliardo seminorm on the right-hand side. It is easy to prove that there exists $c_n>0$ such that 
\begin{equation*}
	\frac{1}{|B|}\int_B |f(x)-f_B|dx\le c_n  r(B)^\delta \left(\frac{1}{|B|}\int_B \int_B \frac{|f(x)-f(y)|^p}{|x-y|^{n+\delta p }}dx dy\right)^\frac{1}{p},
\end{equation*}
for all $p\ge 1$ and $\delta\in (0,1)$. In \cite{BBM1}, J. Bourgain, H. Brezis, and P. Mironescu proved that the previous inequality is far from being optimal. Indeed, the result in \cite{BBM1} yields the following inequality,
\begin{equation}\label{BBM}
	\frac{1}{|B|}\int_B |f(x)-f_B|dx\le c_n (1-\delta)^\frac{1}{p} r(B)^\delta \left(\frac{1}{|B|}\int_B \int_B \frac{|f(x)-f(y)|^p}{|x-y|^{n+\delta p }}dx dy\right)^\frac{1}{p},
\end{equation}
for all ball $B$, $1\le p<\infty$ and $0<\delta<1$. The extra gain factor $(1-\delta)^\frac{1}{p}$ has recently been called the BBM phenomenon \cite{HMPV}. It is a very interesting factor related to the asymptotic behavior of the seminorm when $\delta\to 1^-$ (see \cite{BBM1}),  and it appears naturally in inverse inequalities \cite{Inverse}. We refer also to \cite{MPW} for interesting fractional isoperimetric inequalities with the BBM extra gain factor.

We now present our result on fractional Poincaré-Sobolev inequalities.

\begin{theorem}\label{Fractional}
	Let $n\ge 2$, $\delta\in (0,1)$, $1\le p< \frac{n}{\delta}$ and let $w\in A_1$. Then, there exists a dimensional constant $c_n>0$ such that,
		\begin{align*}
			\left( \frac{1}{w(B)} \int_{B} |f(x)-f_B|^{p \left( \frac{n}{\delta} \right) '} w(x)dx\right)^\frac{1}{p \left( \frac{n}{\delta} \right) '}  \le & c_n \frac{(1-\delta)^\frac{1}{p}}{\delta^{1+\frac{1}{p'}}}  [w]_{A_1}^{\frac{1}{p}} r(B) ^\delta  \\
			 & \cdot  \left( \frac{1}{w(B)} \int_B \int_B \frac{|f(x)-f(y)|^p}{|x-y|^{n+\delta p }} dy \,  w(x)dx \right)^\frac{1}{p}
		\end{align*}
		for any ball $B$.
\end{theorem}

\begin{remark}
  The exponent $\frac{1}{p}$ on the $A_1$ constant is better than the obtained using self-improving methods in \cite{HMPV}. In the case $p=1$, we also improve the quantitative dependence of the $A_1$ constant proved in Theorem 5.7 from \cite{MPW}.  In addition, for weights with large $A_1$ constant, our range of local integrability is larger than the obtained in \cite{HMPV}.  More precisely, $p \left( \frac{n}{\delta} \right) '> \frac{np(1+\log [w]_{A_1})}{n(1+\log [w]_{A_1})-\delta p} $ if and only if $[w]_{A_1}>e^{p-1}$. 
\end{remark}

\subsection{High order derivatives}\label{Section PS m}

For a fixed ball $B$ and an integer $m\in \N$, we consider the space $\mathcal{P}_{m-1}(B)$ of polynomials of degree at most $m-1$ in $n$ variables restricted to the ball $B$.  We can endow the space $\mathcal{P}_{m-1}(B)$ with the following inner product
	\begin{equation*}
		\langle f, g \rangle _B = \frac{1}{|B|}\int_B f(x)g(x)dx.
	\end{equation*}
	
	Let $\{\phi_r \}_r$ be a basis of $\mathcal{P}_{m-1}(B)$ orthonormal with respect to this inner product. We denote by $P_B^{m-1} f$ the projection of a function $f$ onto the space $\mathcal{P}_{m-1}(B)$, 
	\begin{equation*}
		P_B^{m-1} f(x) := \sum_r \left( \frac{1}{|B|} \int_B f(y) \phi_r (y)dy \right) \phi_r (x).
	\end{equation*}
	
	The projection has the following optimality property:
	\begin{equation*}
		\inf _{\pi \in \mathcal{P}_{m-1}(B) } \left( \frac{1}{|B|} \int_B |f(x)-\pi(x)|^p  dx \right)^\frac{1}{p} \approx \left( \frac{1}{|B|} \int_B |f(x)-P_B^{m-1}f(x)|^p  dx \right)^\frac{1}{p} .
	\end{equation*}

	Another important tool that we need is the following subrepresentation formula for high-order derivatives (we refer to \cite{BH93} and also to the book \cite{AH}).

\begin{theorem}\label{formula BH}
	Let $m\ge 1$. There exists a constant $C=C(n,m)>0$ such that for any ball $B$ and any $f\in W^{1,m}(B)$ we have
	\begin{equation*}
		|f(x)-P_B^{m-1}f(x)|\le C I_m (|\nabla^m f| \chi_B)(x),
	\end{equation*}
	for almost every $x\in B$, where $|\nabla^m f| = \sum_{|\alpha|=m} |D^\alpha f|$. 
\end{theorem}

Using this pointwise estimate and Theorem \ref{Teorema I alpha debil 1}, we obtain the following weak Poincaré-Sobolev inequality.

\begin{theorem}
Let $n\ge 2$, $m\ge 1$ and let $1\le p<  \frac{n}{m}$.  Let $w\in A_r$ for some $1\le r \le p$  and consider $q_r$ defined by the relation
	\begin{equation}\label{exp22}
		\frac{1}{p}- \frac{1}{q_r} = \frac{m}{n} \frac{\tau_n [\sigma]_{A_\infty} }{1+r( \tau_n [ \sigma ] _{A_\infty} -1)}
	\end{equation}
	where $\sigma = w^{1-r'}$ and $\tau_n=2^{n+1}$ if $r>1$ and $q_1= p^*_m=\frac{n p }{n- mp }$. Then, there exists a constant $C=C(n, m)>0$ such that,
		\begin{equation*}
			\norm{f-P_B^{m-1}f}_{L^{q_r,\infty}\left(B, \frac{w(x)dx}{w(B)} \right)}  \le C p^*_m [w]_{A_r}^{\frac{1}{p}}  r(B)^m  \left( \frac{1}{w(B)} \int_{B} | \nabla^m f(x)|^p w(x)dx\right)^\frac{1}{p}
		\end{equation*}
		for any ball $B$.
\end{theorem}

The truncation argument does not work for high-order derivatives (see  \cite[p. 121]{H}). In order to obtain strong weighted norm inequalities, we have to use Theorem \ref{coro1}, although it introduces the constant $[\sigma]_{A_\infty}^{1/q_r}$ in the estimates.

\begin{theorem}
	Let $n\ge 2$, $m\ge 1$ and let $1  < p<  \frac{n}{m}$.  Let $w\in A_r$ for some $1\le r \le p$  and consider  $q_r$ defined in \eqref{exp22} if $r>1$ and $q_1=p^*_m= \frac{n p }{n- mp }$. Then, there exists a constant $C=C(n, m)>0$ such that, if $r>1$ we have 
		\begin{equation*}
			\left( \frac{1}{w(B)} \int_{B} |f-P_B^{m-1}f|^{q_r} w\right)^\frac{1}{q_r}  \le C p^*_m [w]_{A_r}^{\frac{1}{p}} [\sigma]_{A_\infty}^\frac{1}{q_r} r(B)^m  \left( \frac{1}{w(B)} \int_{B} | \nabla^m f|^p w\right)^\frac{1}{p}
		\end{equation*}
		for any ball $B$, and if $r=1$ we have 
		\begin{equation*}
			\left( \frac{1}{w(B)} \int_{B} |f-P_B^{m-1}f|^{q_1} w\right)^\frac{1}{q_1}  \le C p^*_m [w]_{A_1}^{\frac{1}{p}}  r(B)^m  \left( \frac{1}{w(B)} \int_{B} | \nabla^m f|^p w\right)^\frac{1}{p}
		\end{equation*}
		for any ball $B$.
\end{theorem}

\subsection{Optimality of the weighted Sobolev exponent}\label{Section exp}

In this section, we address the question of finding the optimal Sobolev exponent for a given weight $w\in A_p$.  In order to state the main result of this section, we need the following definitions. 

\begin{definition}
		Let $w\in A_\infty$, we define 
		\begin{equation*}
			I(w):=\{p\ge 1 :  w\in A_p\}=\left\{
			\begin{array}{ll}
				[1,\infty) & \text{if } w\in A_1 \\
				(\ell_w,\infty) & \text{if } w\in A_\infty \setminus A_1
			\end{array}
		\right.,
		\end{equation*} 
		where 
	\begin{equation*}
		\ell_w : = \inf I(w).
	\end{equation*}
	\end{definition}

We note here that if $w\in A_\infty$ then $w\in A_p$ for each $p>e^{c_n [w]_{A_\infty}}$ (see Theorem 1.3 from \cite{HParissis}). Therefore, for each weight $w\in A_\infty$, the quantity $e^{c_n [w]_{A_\infty}}$ is an upper bound of $\ell_w$. 

\begin{remark}
	Let $w\in A_\infty$, there is a significant difference between $w\in A_1$ and $w\in A_\infty\setminus A_1$. In the first case, $w \in A_{\ell_w}=A_1$. On the other hand, if $w \in A_\infty \setminus A_1$, then $w \notin A_{\ell_w}$. In other words, the infimum $\ell_w$ is a minimum if and only if $w\in A_1$.
\end{remark}

\begin{definition}\label{def p*w}
	Let $1\le p<n \ell_w $ where $w\in A_p$. We define the weighted Sobolev exponent $p^*_w$ by the relation 
\begin{equation*}
	\frac{1}{p}-\frac{1}{p^*_w} = \frac{1}{n}\frac{1}{\ell_w}.
\end{equation*}
\end{definition}

\begin{remark}
	If $w\in A_1$, then $p_w^*=p^*$. 
\end{remark}

We can prove the following theorem using the local bounds for $I_1$ proved throughout the paper. 

\begin{theorem}\label{thm exp 0}
	Let $1\le  p<n \ell_w $ and let $w\in A_p$. Then, for each $1\le q <p^*_w$, we have
	\begin{equation}\label{ps2}
	\left( \frac{1}{w(B)} \int_{B} |f(x)-f_B|^{q} w(x)dx\right)^\frac{1}{q}  \le C_w r(B)  \left( \frac{1}{w(B)} \int_{B} | \nabla f(x)|^p w(x)dx\right)^\frac{1}{p}
\end{equation}
for any ball $B$, where $C_w$ is of the form $c_{n}p^* [w]_{A_r}^\frac{1}{p}$ for some $r>\ell_w$ depending on $q$. 
\end{theorem}

The exponent $p^*_w$ and this type of inequality with the range $p\le q <p^*_w$ are implicit in the proof weighted Poincaré-Sobolev inequalities of \cite{FKS}, as it is remarked in \cite{CMR}. The following result shows that the exponent $p^*_w$ cannot be improved.

\begin{theorem}\label{Best exponent 0}
	We can find a weight $w\in A_p$ such that if the inequality
	\begin{equation*}
	\left( \frac{1}{w(B)} \int_{B} |f(x)-f_B|^{q} w(x)dx\right)^\frac{1}{q}  \le C_w r(B)  \left( \frac{1}{w(B)} \int_{B} | \nabla f(x)|^p w(x)dx\right)^\frac{1}{p}
\end{equation*}
holds for all Lipschitz function $f$, then $q \le p^*_w$.
\end{theorem}

The previous result shows that the inequality \eqref{ps2} is false in general if $q>p^*_w$. On the other hand, we know that it is true in the case $q=p^*_w$ when $w\in A_1\subset A_p$ (Theorem \ref{thm ps 1} with $r=1$). We do not know if it is true at the endpoint $q=p^*_w$ for weights $w\in A_p\setminus A_1$.

We conclude this section by comparing the new optimal exponent $p^*_w$ with the weighted Sobolev exponent $q_r$ defined in \eqref{qr PS eq}. For a weight $w \in A_r$, we have the relation $q_r \le p^*_w$, since $q_r$ is of the form
\begin{equation*}
	\frac{1}{p}-\frac{1}{q_r}= \frac{1}{n} \frac{1}{r-\varepsilon}, 
\end{equation*}
where $\varepsilon>0 $ is given by the left-openness property of the class $A_r$ (see Corollary \ref{precise open property}). By the definition of $\ell_w$ we have $\ell_w < r - \varepsilon$, and then we deduce the relation $q_r < p^*_w$ strictly unless $r = 1$, where $q_1=p^*_w$. The gap between $q_r$ and  $p^*_w$ is given by the precision of the left-openness property: larger $\varepsilon$ (i.e. $r-\varepsilon$ closer to $\ell_w$) narrows the gap, while smaller $\varepsilon$ widens it. Thus, $q_r$ is an approximate optimal exponent given by the left-openness property, whereas $p^*_w$ is the theoretical optimal exponent.

\section{Preliminaries and known results}\label{Section 2}

A weight $w$, a non-negative locally integrable function, is in the Muckenhoupt class $A_p$ with $p>1$ if the quantity
\begin{equation*}
	[w]_{A_p}= \sup_B \left( \frac{1}{|B|}\int_B w(x)dx\right) \left( \frac{1}{|B|}\int_B w(x)^{1-p'}dx\right)^{p-1}
\end{equation*}
is finite, where the supremum is taken over all balls $B$, and $p'$ is the conjugate exponent of $p$, given by $\frac{1}{p}+\frac{1}{p'}=1$. The weight $w$ is in $A_1$ if the quantity
\begin{equation*}
	[w]_{A_1} = \norm{\frac{Mw}{w}}_{L^\infty}
\end{equation*}
is finite, where $M$ is the Hardy-Littlewood maximal operator, 
\begin{equation*}
	Mf(x)=\sup_B \left(\frac{1}{|B|}\int_B |f(y)|dy\right)\chi_B(x).
\end{equation*}

The Hardy-Littlewood maximal operator majorizes several important operators in harmonic analysis. In \cite{M}, Muckenhoupt studied the boundedness of $M$ in weighted $L^p$-spaces, obtaining that for $1\le p<\infty$,
 \begin{equation*}
 	M: L^p(w)\longrightarrow L^{p,\infty}(w)
 \end{equation*}
if and only if $w\in A_p$, and for $1<p<\infty$, 
\begin{equation*}
	M: L^p(w)\longrightarrow L^{p}(w)
\end{equation*}
if and only if $w\in A_p$. The quantities $[w]_{A_p}$ are called $A_p$ constants, and these bounds can be quantified in terms of these constants. We have 
\begin{equation}\label{weak pp Maximal}
	 \norm{M}_{L^p(w)\longrightarrow L^{p,\infty}(w)} \approx_n  [w]_{A_p}^\frac{1}{p}
\end{equation}
for $p\ge 1$, and 
\begin{equation}\label{Buckley}
	\norm{M}_{L^p(w)\longrightarrow L^{p}(w)} \le c_n p' [w]_{A_p}^\frac{1}{p-1}
\end{equation}
for $p>1$. In \cite{B}, Buckley proved that the previous exponents of the $A_p$ constants are sharp. Using the Marcinkiewicz interpolation theorem (see Theorem 1.3.2. from \cite{Grafakos} for example) with the trivial $L^\infty$ strong estimate and the weak type inequality \eqref{weak pp Maximal}, we obtain the following lemma. 
\begin{lemma}\label{lemma maximal interpolada}
	Let $1\le r<p<\infty$, and let  $w\in A_r$, then there exist a dimensional constant $c_n>0$ such that
	\begin{equation}\label{Maximal Lp}
		\| Mf\|_{L^{p}(w)} \le c_n \left(\frac{p}{p-r}\right)^\frac{1}{p} [w]_{A_r}^\frac{1}{p} \|f\|_{L^p(w)}.
	\end{equation}
\end{lemma}

The $A_\infty$ class is defined as the union of all $A_p$ classes, $A_\infty= \bigcup_{p>1}A_p$, this is quantified by the Fuji-Wilson $A_\infty$ constant   
\begin{equation*}
	[w]_{A_\infty} = \sup_B \frac{1}{w(B)} \int_Q M(w\chi_B)(x)dx.
\end{equation*}

In \cite{HPR}, the authors prove that the quantitative bound \eqref{Buckley} can be improved in terms of a mixed $A_p-A_\infty$ constant.

\begin{theorem}[\cite{HPR}]
	Let $1<p<\infty$ and let $w\in A_p$ and $\sigma = w^{1-p'}$. Then there is a dimensional constant $c_n>0$ such that 
	\begin{equation*}\label{mixed}
		\left\| M\right\|_{L^p(w)\longrightarrow L^p(w)} \le c_n p' \left( [w]_{A_p} [\sigma]_{A_\infty} \right) ^\frac{1}{p}.
	\end{equation*}
\end{theorem}

We remark the property $[w]_{A_\infty}\le c_n [w]_{A_p}$ for all $p\ge 1$. We are interested in the left-openness property of $A_p$ weights, which is equivalent to the reverse Hölder inequality. We state the sharp reverse Hölder inequality for $A_\infty$ weights in the Euclidean space obtained in \cite{HPR}.

\begin{theorem}[\cite{HPR}]\label{SRHI}
	Let $w\in A_\infty$ and let $Q$ be a cube. Then
	\begin{equation*}
		\frac{1}{|Q|}\int_Q w(x)^{1+\varepsilon}dx \le 2 \left( \frac{1}{|Q|}\int_Q w(x)dx \right) ^{1+\varepsilon},
	\end{equation*}
	for any $\varepsilon>0$ such that $0<\varepsilon\le \frac{1}{2^{n+1}[w]_{A_\infty}-1}$. 
\end{theorem}

 The range $0<\varepsilon\le \frac{1}{2^{n+1}[w]_{A_\infty}-1}$ is better than the one obtained in \cite{HP} and \cite{HPR} in the context of spaces of homogeneous type. The better range produces an improved precise open property that has been missed in the literature.  

\begin{corollary}\label{precise open property}
	Let $1<p<\infty$ and let $w\in A_p$. We define the quantity 
	\begin{equation*}
		\varepsilon = \frac{p-1}{\tau_n [ \sigma ] _{A_\infty} }, 
	\end{equation*}
	where $\tau_n=2^{n+1}$ and $\sigma = w^{1-p'}$. Then, $w\in A_{p-\varepsilon}$ and 
	\begin{equation*}
		[w]_{A_{p-\varepsilon}} \le 2^{p-1} [w]_{A_p}.
	\end{equation*}
\end{corollary}

\begin{proof}
	Since $w\in A_p$, then $\sigma\in A_{p'}\subset A_\infty$. Let $r(\sigma)=1+\frac{1}{2^{n+1}[\sigma]_{A_\infty}-1}$.  Then, using Theorem \ref{SRHI}, we have
	\begin{equation*}
		\left(\frac{1}{|Q|}\int_Q w(x)dx \right) \left( \frac{1}{|Q|}\int_Q \sigma(x)^{r(\sigma)} dx \right)^\frac{p-1}{r(\sigma)} \le \left(\frac{1}{|Q|}\int_Q w(x)dx \right) \left( \frac{2}{|Q|}\int_Q \sigma(x) dx \right)^{p-1}.
	\end{equation*}
	We may choose $\varepsilon>0$ such that $\frac{p-1}{r(\sigma)}=p-\varepsilon-1$; namely, 
	\begin{equation*}
		\varepsilon = \frac{p-1}{r(\sigma)'} = \frac{p-1}{(1+\frac{1}{2^{n+1}[\sigma]_{A_\infty}-1})'} = \frac{p-1}{2^{n+1}[\sigma]_{A_\infty}} .
	\end{equation*}
	Note that $\varepsilon>0$ and $p-\varepsilon>0$. This yields that $w\in A_{p-\varepsilon}$. 
\end{proof}

We will use the following notation for the weighted Lorentz average over a ball $B$,
\begin{equation*}
	\left\| f\right\|_{L^{p,q}\left(B, \frac{ w(x)dx}{w(B)} \right)} = \left( p  \int_0^\infty t^q \left( \frac{w(\{ x\in B : |f(x)|>t\})}{w(B)}\right)^\frac{q}{p} \frac{dt}{t} \right)^\frac{1}{q}
\end{equation*}
whenever $q<\infty$, and 
\begin{equation*}
	\left\| f\right\|_{L^{p,\infty}\left(B, \frac{ w(x)dx}{w(B)} \right)} = \sup_{t>0} t \left( \frac{w(\{ x\in B : |f(x)|>t\})}{w(B)}\right)^\frac{1}{p} .
\end{equation*}

We also need the following lemma from \cite{HMPV}. The result can be found on \cite{FLW} in the case $r=1$ and \cite{FH} when $r>1$, but we need the version of \cite{HMPV} where the precise dependence of the parameters is considered. 

\begin{lemma}\label{rep}
	Let $B_0$ be a ball in $\R^n$. Assume that $0<\alpha <n $ and consider $0<\eta <n-\alpha$ and $1 \le r <\infty$. Suppose that there exists a constant $\kappa>0$, a function $f\in L^1(B_0)$ and a non-negative measurable function $g$ such that 
	\begin{equation*}
		\frac{1}{|B|}\int_B |f(x)-f_B|dx \le \kappa r(B)^\alpha \left(\frac{1}{|B|}\int_B g(x)^r dx \right)^\frac{1}{p}
	\end{equation*}
	for every ball $B\subset B_0$. Then, there exists a dimensional constant $C_n$ such that 
	\begin{equation*}
		|f(x)-f_{B_0}|\le C_n \frac{\kappa}{\alpha^\frac{1}{r'} \eta^\frac{1}{r}} r(B_0)^\frac{\alpha}{r'} \left( I_\alpha (g^r \chi_{B_0})(x)\right)^\frac{1}{r}
	\end{equation*} 
	for every Lebesgue point $x\in B_0$ of $f$. 
\end{lemma}

As usual, $C$ denotes various positive constants, and $C(\alpha, \beta,...)$ or $C_{\alpha, \beta, ...}$, denote such constants depending only on $\alpha, \beta,...$ . These constants can change even in the same string of estimates.

\section{Weighted Hedberg's inequality and consequences}\label{Section 3}

We need the following weighted Hedberg's inequality (see \cite{Hedberg} for the classical case).

\begin{lemma}\label{lemma Hedberg}
	Let $n\ge 1$, $\alpha \in (0,n)$, let $1< p<\frac{n}{\alpha}$ and let $w$ be a weight. Let $1< r \le  p $ and consider $q$ defined by the relation
	\begin{equation*}
		\frac{1}{p}- \frac{1}{q} = \frac{\alpha}{n} \frac{1}{r}.
	\end{equation*}
	Then, there exists a dimensional constant $c_n>0$ such that for any ball $B$ such that $0<\int_B w(y) ^{1-r'} dy <\infty$ and every function $f\in L^p(B, w)$, we have 
	\begin{equation*} 
		|I_\alpha f(x)| \le \frac{c_n}{\alpha } p^*_\alpha   Mf(x) ^\frac{p}{q}  \left\| f \right\| _{L^p(B, w)}^{1-\frac{p}{q}} \left( \int_B w(y) ^{1-r'} dy \right) ^\frac{\alpha }{nr'},
	\end{equation*}
	for any $x\in B$. 
\end{lemma}

This inequality was first proved in \cite{CF} in the case $\alpha =1$ and $w\in A_p$. The proof of this inequality for $w\in A_p$ is essentially contained in the proof of Theorem 1.1 from \cite{ACS}; we have chosen to include the proof because we need to be precise with the main parameters involved.

\begin{proof}
	Let $B$ be a ball such that $0<\int_B w(y) ^{1-r'} dy <\infty$, and let $f\in L^p(B, w)$. Fix $x\in B$ and let $\delta>0$. We can decompose as follows,
	\begin{equation*}
		|I_\alpha f(x)| \le \int_{ |x-y|\le \delta} \frac{|f(y)|}{|x-y|^{n-\alpha}}dy  + \int_{ |x-y|> \delta} \frac{|f(y)|}{|x-y|^{n-\alpha}}dy = I+II. 
	\end{equation*}
	
	We will estimate each term separately. Firstly,
	\begin{align*}
		I= & \sum_{k=0}^\infty \int_{ \frac{\delta}{2^{k+1}} < |x-y|\le \frac{\delta}{2^k}}  \frac{|f(y)|}{	|x-y|^{n-\alpha}}dy  \\
		\le & \sum_{k=0}^\infty \left( \frac{2^{k+1}}{\delta}\right) ^{n-\alpha} \int_{  |x-y|\le \frac{\delta}{2^k}} |f(y)| dy  \\
		\le & 2^{n-\alpha} \omega_n \delta^\alpha Mf(x)  \sum_{k=0}^\infty \frac{1}{2^{\alpha k }}  \\ 
		= &  \frac{2^n \omega_n}{2^\alpha-1} \delta^\alpha Mf(x),
	\end{align*}
	where $\omega_n$ is the Lebesgue measure of the unit ball. Secondly, using Hölder's inequality, we obtain
	\begin{align*}
		II= & \int_{ |x-y|> \delta} \frac{|f(y)|}{|x-y|^{n-\alpha}} w(y)^\frac{1}{p}w(y)^{-\frac{1}{p}}  dy  \\
		\le & \norm{f}_{L^p(B, w)} \left( \int_{B \cap \{ |x-y|> \delta\}} \frac{w(y)^{1-p'} }{|x-y|^{(n-\alpha)p'}}  dy\right) ^\frac{1}{p'}.
	\end{align*}

	Let $1<r<p$, we use Hölder's inequality again with $\frac{1-r'}{1-p'}>1$,
	\begin{align*}
		\int_{B \cap \{ |x-y|> \delta\}} \frac{w(y)^{1-p'}  }{|x-y|^{(n-\alpha)p'}} dy \le & \left( \int_{B } w(y)^{1-r'}  dy  \right)^\frac{1-p'}{1-r'}  \\
		& \cdot \left( \int_{B \cap \{ |x-y|> \delta\}} \frac{dy }{|x-y|^{(n-\alpha)p' \left( \frac{1-r'}{1-p'}\right)'}}   \right)^\frac{1}{\left( \frac{1-r'}{1-p'}\right)'}	\\
		= & \left( \int_{B } w(y)^{1-r'}  dy  \right)^\frac{r-1}{p-1}  \left( \int_{B \cap \{ |x-y|> \delta\}} \frac{dy}{|x-y|^{ \frac{p(n-\alpha)}{p-r} }}    \right)^\frac{p-r}{p-1}.
	\end{align*}
	
	Hence, 
	\begin{equation*}
		II\le \norm{f}_{L^p(B, w)} \left( \int_{B } w(y)^{1-r'}  dy  \right)^\frac{r-1}{p} \left( \int_{B \cap \{ |x-y|> \delta\}} \frac{dy}{|x-y|^{ \frac{p(n-\alpha)}{p-r} }}    \right)^\frac{p-r}{p}.
	\end{equation*}
	
	We can estimate the second integral,
	\begin{align*}
		\left( \int_{B \cap \{ |x-y|> \delta\}} \frac{dy}{|x-y|^{ \frac{p(n-\alpha)}{p-r} }}   \right)^\frac{p-r}{p} = & \left( \int_0^\infty \left| \left\lbrace y\in B : \delta < |x-y|< t^{- \frac{p-r}{p(n-\alpha)} } \right\rbrace \right| dt \right)^\frac{p-r}{p}\\
		= & \left( \int_0^{\delta^{ - \frac{p(n-\alpha) }{p-r} }} \left| \left\lbrace y\in B : \delta < |x-y|< t^{- \frac{p-r}{p(n-\alpha) }} \right\rbrace \right| dt\right)^\frac{p-r}{p} \\
		\le  & \left(\int_0^{\delta^{ - \frac{p(n-\alpha) }{p-r} }} \left| B \left(x, t^{- \frac{p-r}{p(n-\alpha) }} \right) \right| dt \right)^\frac{p-r}{p} \\
		= & \left(\omega_n \int_0^{\delta^{ - \frac{p(n-\alpha) }{p-r} }} t^{- \frac{(p-r)n}{p(n-\alpha) }}  dt \right)^\frac{p-r}{p} \\
		\le & c_n \frac{p(n-\alpha)}{n-\alpha p}  \delta^{\alpha - \frac{nr}{p}}.  
	\end{align*}
	
	If $r=p$, we obtain a similar estimate,
	\begin{equation*}
		II\le \norm{f}_{L^p(B, w)} \delta^{\alpha-n} \left( \int_{B } w(y)^{1-p'}   dy\right) ^\frac{1}{p'}.
	\end{equation*}
	
	Therefore, 
	\begin{equation*}
		|I_\alpha f(x)|  \le  C \left( \delta^\alpha Mf(x) + \delta^{\alpha - \frac{nr}{p}} \norm{f}_{L^p(B, w)} \left( \int_{B } w(y)^{1-r'}  dy  \right)^\frac{r-1}{p} \right),
	\end{equation*}
	where $C= c_n \max \left\lbrace  \frac{1}{2^\alpha-1}  ,  \frac{p(n-\alpha)}{n-\alpha p}\right\rbrace  $. Observe that since $\frac{\log 2}{2^\alpha- 1}\le \frac{1}{ \alpha   }$, we have 
	\begin{equation*}
		C \le c_n p \left(\frac{1}{\alpha }+ \frac{n-\alpha }{n-\alpha p } \right) \le \frac{c_n p}{\alpha (n-\alpha p)}=\frac{c_n}{\alpha} p^*_\alpha . 
	\end{equation*}
	Without loss of generality, we may assume that $Mf(x)>0$, since otherwise $f=0$ and there is nothing to prove. The choice 
	\begin{equation*}
		\delta = \left( \frac{\norm{f}_{L^p(B, w)} \left( \int_{B } w(y)^{1-r'}  dy  \right)^\frac{r-1}{p}}{Mf(x)} \right) ^\frac{p}{nr}
	\end{equation*}
	yields the desired inequality 
	\begin{equation*} 
		|I_\alpha f(x)| \le  \frac{c_n}{\alpha} p^*_\alpha   Mf(x) ^\frac{p}{q}  \left\| f \right\| _{L^p(B, w)}^{1-\frac{p}{q}} \left( \int_B w(y) ^{1-r'} dy \right) ^\frac{\alpha }{nr'}
	\end{equation*}
	for each $x\in B$. 
\end{proof}

We will prove the following result, which implies Theorem \ref{coro1}.

\begin{proposition}\label{Teorema I alpha fuerte}
	Let $n \ge 1$, $\alpha \in (0,n)$ and $1<p<\frac{n}{\alpha}$. Let $w\in A_p$ and consider $q_r$ defined by the relation
	\begin{equation}\label{exp33}
		\frac{1}{p}- \frac{1}{q_r} = \frac{\alpha}{n} \frac{1}{ r} 
	\end{equation}
	where $1\le r<p$. If $r\in I(w)$, then there exists a dimensional constant $c_n>0$ such that  
	\begin{equation}\label{eq42}
		\left( \frac{1}{w(B)} \int_{B} |I_\alpha f|^{q_r} w\right)^\frac{1}{q_r} \le \frac{c_n}{\alpha} p^*_\alpha  \left(\frac{p}{p-r}\right)^\frac{1}{q_r} [w]_{A_r}^{\frac{1}{p}}  r(B)^\alpha  \left( \frac{1}{w(B)} \int_{B} | f|^p w\right)^\frac{1}{p},
	\end{equation}
	for any ball $B$ and every function $f\in L^p(B, w)$. 
\end{proposition}

\begin{proof}[Proof]
	Let $r\in I(w)$. Let assume first that $1<r<p$, then we have $0<\int_B w(y) ^{1-r'} dy <\infty$, since $w\in A_r$. Applying Lemma \ref{lemma Hedberg} we obtain
	\begin{equation*}
		|I_\alpha f(x)| \le \frac{c_n}{\alpha} p^*_\alpha   Mf(x) ^\frac{p}{q_r}  \left\| f \right\| _{L^p(B, w)}^{1-\frac{p}{q_r}} \left( \int_B w(y) ^{1-r'} dy \right) ^\frac{\alpha }{nr'},
	\end{equation*}
	for all $x\in B$. Therefore,
	\begin{align*}
		\left( \frac{1}{w(B)} \int_{B} |I_\alpha f|^{q_r} w \right)^\frac{1}{q_r} \le & \frac{c_n}{\alpha} p^*_\alpha  \left\| f \right\| _{L^p(B, w)}^{1-\frac{p}{q_r}} \left( \int_B w ^{1-r'} \right) ^\frac{\alpha }{nr'} \left( \frac{1}{w(B)} \int_{B} |Mf |^p w \right)^\frac{1}{q_r}\\
		= & \frac{c_n}{\alpha} p^*_\alpha  \left\| f \right\| _{L^p(B, w)}^{1-\frac{p}{q_r}} \frac{w(B)^\frac{\alpha}{nr}}{w(B)^\frac{1}{p}}  \left( \int_B w ^{1-r'} \right) ^\frac{\alpha (r-1) }{nr} \left( \int_{B} |Mf |^p w \right)^\frac{1}{q_r}\\
		\le & \frac{c_n}{\alpha} p^*_\alpha  [w]_{A_r}^\frac{\alpha }{nr} r(B)^\alpha \frac{1}{w(B)^\frac{1}{p}} \left\| f \right\| _{L^p(B, w)}^{1-\frac{p}{q_r}}  \left( \int_{\R^n} |Mf |^p w \right)^\frac{1}{q_r}
	\end{align*}
	where we have used the $A_r$ condition. Since $w\in A_r$ with $ r<p<\infty$, we  use the estimate \eqref{Maximal Lp},
	\begin{align*}
		\left( \frac{1}{w(B)} \int_{B} |I_\alpha f|^{q_r} w \right)^\frac{1}{q_r} \le &  \frac{c_n}{\alpha} p^*_\alpha  \left(\frac{p}{p-r}\right)^\frac{1}{q_r} [w]_{A_r}^\frac{\alpha }{nr} r(B)^\alpha \frac{1}{w(B)^\frac{1}{p}} \left\| f \right\| _{L^p(B, w)}^{1-\frac{p}{q_r}}  [w]_{A_r}^{ \frac{1}{q_r}} \left\| f \right\| _{L^p(B, w)}^\frac{p}{q_r}\\
		= & \frac{c_n}{\alpha} p^*_\alpha  \left(\frac{p}{p-r}\right)^\frac{1}{q_r} [w]_{A_r}^{\frac{1}{p}}  r(B)^\alpha \left( \frac{1}{w(B)} \int_{B} | f|^p w\right)^\frac{1}{p}.
	\end{align*}
	This yields \eqref{eq42}. Let us consider the remaining case $r=1\in I(w)$. We apply the classical Hedberg's inequality (see \cite{Hedberg}) and we use $w\in A_1$,       
	\begin{align*}
		|I_\alpha f(x)| \le & \frac{c_n}{\alpha} p^*_\alpha   Mf(x) ^\frac{p}{q_1}  \left\| f \right\| _{L^p(B, dx)}^{1-\frac{p}{q_1}} \le  \frac{c_n}{\alpha} p^*_\alpha  [w]_{A_1}^\frac{\alpha}{n} r(B)^\alpha   Mf(x) ^\frac{p}{q_1}  \left\| f \right\| _{L^p\left( B, \frac{w(x)dx}{w(B)} \right) }^{1-\frac{p}{q_1}}  .
	\end{align*}
	Using again \eqref{Maximal Lp} with $r=1$, we have,
	\begin{align*}
		\| I_\alpha f\|_{L^{q_1}(B, w)} \le & \frac{c_n}{\alpha} p^*_\alpha [w]_{A_1}^\frac{\alpha}{n} r(B)^\alpha  \left\| f \right\| _{L^p\left( B, \frac{w(x)dx}{w(B)} \right) }^{1-\frac{p}{q_1}}  \| Mf\|_{L^{p}(B, w)} ^\frac{p}{q_1}\\
		\le & \frac{c_n}{\alpha} p^*_\alpha \left(p'\right)^\frac{1}{q_1} [w]_{A_1}^\frac{\alpha}{n} r(B)^\alpha  \left\| f \right\| _{L^p\left( B, \frac{w(x)dx}{w(B)} \right) }^{1-\frac{p}{q_1}} [w]_{A_1}^{\frac{1}{q_1}} \| f\|_{L^{p}(B, w)} ^\frac{p}{q_1}\\
		= & \frac{c_n}{\alpha} p^*_\alpha  \left(p'\right)^\frac{1}{q_1} [w]_{A_1}^\frac{1}{p} r(B)^\alpha w(B)^\frac{1}{q_1}  \left\| f \right\| _{L^p\left( B, \frac{w(x)dx}{w(B)} \right) },
	\end{align*}
	and this yields the desired inequality \eqref{eq42} which completes the proof. 
\end{proof}

\begin{remark}
	It should be noted that a similar inequality holds for pairs of weights $(u,v)\in A_r$ with $1\le r<p$ (we refer to \cite{Journe} for the precise definition). 
The proof follows the same path, using a version of Lemma \ref{lemma maximal interpolada} for pairs of weights $(u,v)\in A_r$. 
\end{remark}

Using the improved precise open property (Corollary \ref{precise open property}) and the previous proposition, we are able to prove Theorem \ref{coro1}.

\begin{proof}[Proof of Theorem \ref{coro1}]
	We may assume $r>1$ since the case $r=1$ has been proved in the previous proposition.  Using Corollary \ref{precise open property}, since $w\in A_r$, then $w\in A_{r-\varepsilon}$ where  
	\begin{equation*}
		\varepsilon= \frac{r-1}{\tau_n [ \sigma ] _{A_\infty} }, 
	\end{equation*}
	and $[w]_{A_{r-\varepsilon}} \le 2^{r-1} [w]_{A_r}$. Observe that 
	\begin{equation*}
		r-\varepsilon = r - \frac{r-1}{\tau_n [ \sigma ] _{A_\infty} } =\frac{1+r(\tau_n [\sigma]_{A_\infty} -1)}{\tau_n [\sigma]_{A_\infty} }.
	\end{equation*} 
	Since $w\in A_{r-\varepsilon}$, then $r-\varepsilon\in I(w)$ and we can apply Proposition \ref{Teorema I alpha fuerte},
	\begin{align*}
		\left( \frac{1}{w(B)} \int_{B} |I_\alpha f|^{q} w\right)^\frac{1}{q} \le & \frac{c_n}{\alpha} p^*_\alpha  \left(\frac{p}{p-(r-\varepsilon)}\right)^\frac{1}{q} [w]_{A_{r-\varepsilon}}^{\frac{1}{p}}  r(B)^\alpha  \left( \frac{1}{w(B)} \int_{B} | f|^p w\right)^\frac{1}{p}\\
		\le & \frac{c_n}{\alpha} p^*_\alpha 2^{\frac{1}{p'}}  \left(\frac{p}{p-r+\varepsilon}\right)^\frac{1}{q} [w]_{A_{r}}^{\frac{1}{p}}  r(B)^\alpha  \left( \frac{1}{w(B)} \int_{B} | f|^p w\right)^\frac{1}{p}\\
		\le & \frac{c_n}{\alpha} p^*_\alpha \left(\frac{p}{p-1} \tau_n [\sigma]_{A_\infty} \right)^\frac{1}{q} [w]_{A_{r}}^{\frac{1}{p}}  r(B)^\alpha  \left( \frac{1}{w(B)} \int_{B} | f|^p w\right)^\frac{1}{p}\\
		= & \frac{c_n}{\alpha} p^*_\alpha (p')^\frac{1}{q}     [w]_{A_{r}}^{\frac{1}{p}} [\sigma]_{A_\infty} ^\frac{1}{q} r(B)^\alpha  \left( \frac{1}{w(B)} \int_{B} | f|^p w\right)^\frac{1}{p},
	\end{align*}
	 where we have used the bound
	\begin{align*}
		\frac{p}{p-r+\varepsilon} = & \frac{p \tau_n [\sigma]_{A_\infty} }{p \tau_n [\sigma]_{A_\infty} -1 - r(\tau_n [\sigma]_{A_\infty} -1)} \\\le & \frac{p \tau_n [\sigma]_{A_\infty} }{p \tau_n [\sigma]_{A_\infty} -1 - p(\tau_n [\sigma]_{A_\infty} -1)} \\
		= & \frac{p}{p-1} \tau_n [\sigma]_{A_\infty},
	\end{align*}
	since $r\le p$ and $\tau_n [\sigma]_{A_\infty} >1$. This completes the proof. 
\end{proof}

As before, we need the following proposition, which implies Theorem \ref{Teorema I alpha debil 1}.

\begin{proposition}\label{Teorema I alpha debil}
	Let $n\ge 1$, $\alpha \in (0,n)$ and $1\le p<\frac{n}{\alpha}$. Let $w\in A_p$ and consider $q_r$ defined in \eqref{exp33} with $1\le r \le p$. If $r\in I(w)$, then there exists a dimensional constant $c_n>0$ such that  
	\begin{equation*}
		\left\|I_\alpha f\right\|_{ L^{q_r,\infty} \left( B , \frac{w(x)dx}{w(B)} \right) } \le \frac{c_n}{\alpha} p^*_\alpha [w]_{A_r}^{\frac{1}{p}}  r(B)^\alpha  \left( \frac{1}{w(B)} \int_{B} | f(x)|^p w(x)dx\right)^\frac{1}{p},
	\end{equation*}
	for any ball $B$ and every function $f\in L^p(B, w)$. 
\end{proposition}

The proof of the previous proposition is analogous to the proof of Proposition \ref{Teorema I alpha fuerte}, using \eqref{weak pp Maximal} instead of \eqref{Maximal Lp}, which produces a better constant. We remark that if $p=1$, we can use the classical Hedberg's inequality, and an analogous argument produces the desired inequality for each $w\in A_1$.

\begin{proof}[Proof of Theorem \ref{Teorema I alpha debil 1}]
	We may assume $r>1$ since the case $r=1$ is contained in the previous proposition. Using the left-openness property stated in Corollary \ref{precise open property}, since $w\in A_r$, then $w\in A_{r-\varepsilon}$ where   
		\begin{equation*}
		\varepsilon= \frac{r-1}{\tau_n [ \sigma ] _{A_\infty} }, 
	\end{equation*}
	and we have $[w]_{A_{r-\varepsilon}} \le 2^{r-1} [w]_{A_r}$. Note that $q=q_{r-\varepsilon}$. Since $w\in A_{r-\varepsilon}$, then $r-\varepsilon\in I(w)$ and we can apply Proposition \ref{Teorema I alpha debil},
	\begin{align*}
		\left\|I_\alpha f\right\|_{ L^{q,\infty} \left( B , \frac{w(x)dx}{w(B)} \right) }  \le & \frac{c_n}{\alpha} p^*_\alpha  [w]_{A_{r-\varepsilon}}^{\frac{1}{p}}  r(B)^\alpha  \left( \frac{1}{w(B)} \int_{B} | f(x)|^p w(x)dx\right)^\frac{1}{p}\\
		\le & \frac{c_n}{\alpha} p^*_\alpha  [w]_{A_{r}}^{\frac{1}{p}}  r(B)^\alpha  \left( \frac{1}{w(B)} \int_{B} | f(x)|^p w(x)dx\right)^\frac{1}{p}.
	\end{align*}
\end{proof}

\section{Poincaré-Sobolev inequalities and weighted Sobolev exponent}\label{Section 4}

In this section, we provide the proofs of the results stated in Sections \ref{Section PS}, \ref{Section PS m}, and \ref{Section exp}.

\begin{proof}[Proof of Theorem \ref{thm ps 1}]
	Since $w\in A_r$ we have $r\in I(w)$, and we can apply Theorem \ref{Teorema I alpha debil 1} with $\alpha=1$ to obtain 
	\begin{align*}
		\left\|f-f_B\right\|_{ L^{q_r,\infty} \left( B , \frac{w(x)dx}{w(B)} \right) } \le & c_n \left\|I_1(|\nabla f|\chi_B )\right\|_{ L^{q_r,\infty} \left( B , \frac{w(x)dx}{w(B)} \right) } \\
		\le & c_n  p^* [w]_{A_r}^{\frac{1}{p}}  r(B)  \left( \frac{1}{w(B)} \int_{B} | \nabla f(x)|^p w(x)dx\right)^\frac{1}{p},
	\end{align*} 
	where in the first inequality, we have used \eqref{subrepresentation formula 1}.  We obtain the desired inequalities using the truncation method (see \cite{H, KO03}). 
\end{proof}

The proofs of Poincaré-Sobolev inequalities for high-order derivatives follow the same ideas using Theorem \ref{formula BH} instead of \eqref{subrepresentation formula 1}.

\begin{proof}[Proof of Theorem \ref{Fractional}]
	Fix a ball $B$, using \eqref{BBM} and  Lemma \ref{rep} with $\alpha=\delta$,  $\eta= \frac{n-\delta}{2} $ and $r=p$, we obtain the following representation formula,
	\begin{equation}\label{fractional rep formula}
		|f(x)-f_B|\le c_{n} \frac{(1-\delta)^\frac{1}{p}}{\delta^\frac{1}{p'}}  r(B)^\frac{\delta}{p'} I_\delta (g_{f,B}^p\chi_B)(x)^\frac{1}{p} ,
	\end{equation}
	for almost every $x\in B$, where 
	\begin{equation*}
		g_{f,B}(x)=\left( \int_B  \frac{|f(x)-f(y)|^p}{|x-y|^{n+\delta p }} dy \right)^\frac{1}{p}.
	\end{equation*}
	Using \eqref{fractional rep formula} and Proposition \ref{Teorema I alpha debil}  with $\alpha = \delta$, $p=r=1$ and noting that $q_1= \left( \frac{n}{\delta} \right) '$   we have,
	\begin{align*}
		\norm{ f-f_B }_{L^{p \left( \frac{n}{\delta} \right) ', \infty} \left( B, \frac{w(x)dx}{w(B)}\right) }\le & c_{n} \frac{(1-\delta)^\frac{1}{p}}{\delta^\frac{1}{p'}}  r(B)^\frac{\delta}{p'} \norm{ I_\delta (g_{f,B}^p\chi_B)^\frac{1}{p} }_{L^{p \left( \frac{n}{\delta} \right) ', \infty} \left( B, \frac{w(x)dx}{w(B)}\right) } \\
		= & c_{n} \frac{(1-\delta)^\frac{1}{p}}{\delta^\frac{1}{p'}}  r(B)^\frac{\delta}{p'} \norm{ I_\delta (g_{f,B}^p\chi_B) }_{L^{ \left( \frac{n}{\delta} \right) ', \infty} \left( B, \frac{w(x)dx}{w(B)}\right) }^\frac{1}{p}\\
		\le &  \frac{c_n   }{\delta  (n-\delta )} \frac{(1-\delta)^\frac{1}{p}}{\delta^\frac{1}{p'}}  r(B)^\frac{\delta}{p'} [w]_{A_1}^\frac{1}{p} r(B)^\frac{\delta}{p}  \left( \frac{1}{w(B)}\int_B |g_{f,B}(x)|^p w(x)dx \right)^\frac{1}{p}\\
		\le & c_n  \frac{(1-\delta)^\frac{1}{p}}{\delta^{1+\frac{1}{p'}}}   [w]_{A_1}^\frac{1}{p} r(B)^\delta   \left( \frac{1}{w(B)} \int_B \int_B \frac{|f(x)-f(y)|^p}{|x-y|^{n+\delta p }} dy \,  w(x)dx \right)^\frac{1}{p}.
	\end{align*}
	The conclusion of the theorem comes from the truncation method, which is also available in the fractional setting (we refer to \cite{DLV}). 
\end{proof}

\begin{proof}[Proof of Theorem \ref{thm exp 0}]
	
	Let $1\le p < n\ell_w$ where $w\in A_p$.  We may assume without loss of generality that $q_p< q <p^*_w$, where $q_p$ is defined by \eqref{exp33}, because to prove the remaining cases, we use this case and Jensen's inequality. Then, there exists $r\in (\ell_w, p)$ such that $q=q_r$, where
	\begin{equation*}
		\frac{1}{p}-\frac{1}{q_r}=\frac{1}{n}\frac{1}{r}.
	\end{equation*}
	Since $\ell_w < r$, then $w\in A_r$ and using \eqref{subrepresentation formula 1} and Proposition \ref{Teorema I alpha debil} we have 
	\begin{align*}
		\left\|f-f_B\right\|_{ L^{q_r,\infty} \left( B , \frac{w(x)dx}{w(B)} \right) }  \le & c_n \left\|I_1 (|\nabla f|\chi_B)\right\|_{ L^{q_r,\infty} \left( B , \frac{w(x)dx}{w(B)} \right) } \\
		\le & c_n p^* [w]_{A_r}^{\frac{1}{p}}  r(B)  \left( \frac{1}{w(B)} \int_{B} |\nabla  f(x)|^p w(x)dx\right)^\frac{1}{p},
	\end{align*}
	for any ball. Using the truncation argument, we obtain the desired inequality \eqref{ps2}.
\end{proof}

We will prove the following result, which is equivalent to Theorem \ref{Best exponent 0}. 

\begin{proposition}\label{Best exponent}
	We can find a weight $w\in A_p$ such that if the inequality 
	\begin{equation}\label{eq51}
	\left( \frac{1}{w(B)} \int_{B} |f(x)-f_B|^{q} w(x)dx\right)^\frac{1}{q}  \le C_w r(B)  \left( \frac{1}{w(B)} \int_{B} | \nabla f(x)|^p w(x)dx\right)^\frac{1}{p}
\end{equation}
holds for all Lipschitz function $f$, where $q$ is defined by 
	\begin{equation*}
		\frac{1}{p}-\frac{1}{q} = \frac{1}{n(\ell_w-\gamma)}
	\end{equation*}
	with $\gamma$ such that $n(\ell_w-\gamma)>0$. Then $\gamma\le 0$. 
\end{proposition}

We will follow the ideas of the proof of Proposition 7.4 from \cite{CPER}. We need the following lemma from \cite{CPER}.

\begin{lemma}$($\cite[Lemma 7.1]{CPER}$)$
	Let $\mu$ be a finite measure such that $\operatorname{supp} (\mu)\subset \Omega\subset \R^n$. Consider a subset $E\subset \Omega$ such that $\mu(E)\ge \lambda \mu(\Omega)$ for some $\lambda\in (0,1)$ and a function $f$ vanishing on $E$. Then, for any constant $a\in \R$ we have 
	\begin{equation*}
		\left\| a \right\|_{L^q(\mu)} \le \frac{1}{\lambda^q} \left\| f- a \right\| _{L^q(\mu)} .
	\end{equation*} 
\end{lemma}

\begin{proof}[Proof of Theorem \ref{Best exponent}]
	Consider $B=B(0,1)\subset \R^n$ with $n\ge 2$. Let $1<p<n$, and consider the weight $w(x)=|x|^{\delta-n}$ with $n<\delta<np$, then $w\in A_p$ (see Example 7.1.7. in \cite{Grafakos}). Concretely, $\ell_w= \frac{\delta}{n}>1$, and this implies that $w\in A_r$ for all $r\in (\frac{\delta}{n}, \infty)$. Observe that $w(B)= \frac{\omega_{n-1}}{\delta}$, where $\omega_{n-1}$ is the $n-1$ dimensional measure of $\mathbb{S}^{n-1}$. 
 	Let $0<\varepsilon<\frac{1}{2}$ be a small number which we will choose later, and let define $E=B\setminus B(0,2\varepsilon)$. We define on $B$ a piecewise affine Lipschitz function $f$ such that $f(x)=0$ for all $x\in E$ and $f(x)=1$ for all $x\in B(0,\varepsilon)$. Note that $|\nabla f(x)| = \frac{1}{\varepsilon}$ for all $x\in B(0,2\varepsilon)\setminus B(0,\varepsilon)$ and it is 0 otherwise. In order to apply the previous lemma, we compute
	\begin{align*}
		w(E)=  \int_E |x|^{\delta-n} dx = \omega_{n-1} \int _{2\varepsilon}^1 \rho^{\delta+n-n}\frac{d\rho}{\rho } = \frac{\omega_{n-1}}{\delta} \left( 1-(2\varepsilon)^\delta\right) = \left( 1-(2\varepsilon)^\delta\right) w(B). 
	\end{align*}
	Then, applying the lemma and that inequality \eqref{eq51} holds, we have
	\begin{align*}
		\left( \frac{1}{w(B)} \int_{B} |f(x)|^{q} w(x)dx\right)^\frac{1}{q} \le & \left( \frac{1}{w(B)} \int_{B} |f(x)-f_B|^{q} w(x)dx\right)^\frac{1}{q} + \left( \frac{1}{w(B)} \int_{B} |f_B|^{q} w(x)dx\right)^\frac{1}{q} \\
		\le & \left( 1+ \frac{1}{\left(  1-(2\varepsilon)^\delta \right)^q } \right) \left( \frac{1}{w(B)} \int_{B} |f(x)-f_B|^{q} w(x)dx\right)^\frac{1}{q}\\
		\le &   \left( 1+ \frac{1}{\left(  1-(2\varepsilon)^\delta \right)^q } \right)  C_w r(B)  \left( \frac{1}{w(B)} \int_{B} | \nabla f(x)|^p w(x)dx\right)^\frac{1}{p}.
	\end{align*}
	
	We will study each term separately, 
	\begin{align*}
		\left( \frac{1}{w(B)} \int_{B} |f(x)|^{q} w(x)dx\right)^\frac{1}{q} & \ge   \left( \frac{\delta}{\omega_{n-1}} \int_{B(0,\varepsilon)} |f(x)|^q w(x) dx\right)^\frac{1}{q} \\
		& = \left( \frac{\delta}{\omega_{n-1}}  \int_{B(0,\varepsilon)}  |x|^{\delta-n} dx\right)^\frac{1}{q} \\
		& =  \left( \frac{\omega_{n-1} \delta}{\omega_{n-1} \delta}   \varepsilon^{\delta}  \right)^\frac{1}{q} \\
		& =   \, \varepsilon^\frac{\delta}{q}.
	\end{align*}
	
	On the other hand, 
	\begin{align*}
		 r(B)  \left( \frac{1}{w(B)} \int_{B} | \nabla f(x)|^p w(x)dx\right)^\frac{1}{p} & =   \left( \frac{\delta}{\omega_{n-1}} \int_{B(0,2\varepsilon)\setminus B(0,\varepsilon) } \frac{1}{\varepsilon^p}  |x|^{\delta-n} dx\right)^\frac{1}{p}\\
		&\le  \frac{1}{\varepsilon} \left(\frac{\delta}{\omega_{n-1}}  c_n  2^{\delta -1} \varepsilon^\delta  \right)^\frac{1}{p} \\
		&=    2^{(\delta-1)\frac{1}{p}} \delta^\frac{1}{p} \varepsilon^{\frac{\delta}{p}-1}  \\
		& \lesssim_{n}   \delta^\frac{1}{p} \varepsilon^{\frac{\delta}{p}-1}  .
	\end{align*}
	
	Hence,
	\begin{equation}\label{optimality eq1}
		\varepsilon^\frac{\delta}{q} \lesssim_{n} C_w  \left( 1+ \frac{1}{\left(  1-(2\varepsilon)^\delta \right)^q } \right) \delta^\frac{1}{p} \varepsilon^{\frac{\delta}{p}-1}  ,
	\end{equation}	
	and this implies,
	\begin{equation*} 
		1 \lesssim_{n} C_w \left( 1+ \frac{1}{\left(  1-(2\varepsilon)^\delta \right)^q } \right) \varepsilon^{\frac{\delta}{p}-\frac{\delta}{q}-1}  .
	\end{equation*}
	
	Observe that, since $\ell_w=\frac{\delta}{n}$, we have
	\begin{equation*}
		\frac{1}{p}-\frac{1}{q}= \frac{1}{n (\ell_w-\gamma)}=\frac{1}{n \left( \frac{\delta}{n} -\gamma \right)} = \frac{1}{\delta-\gamma n }.
	\end{equation*}
	
	Therefore, 
	\begin{align*}
		1 \le & C_w \left( 1+ \frac{1}{\left(  1-(2\varepsilon)^\delta \right)^q } \right)   \varepsilon^{\frac{\delta}{p}-\frac{\delta}{q}-1} \\
		= & C_w \left( 1+ \frac{1}{\left(  1-(2\varepsilon)^\delta \right)^q } \right)  \varepsilon^{\frac{\delta}{\delta-\gamma n}-1}\\
		= & C_w \left( 1+ \frac{1}{\left(  1-(2\varepsilon)^\delta \right)^q } \right)  \varepsilon^{\frac{\gamma n}{\delta-\gamma n}} .
	\end{align*}
	
	Since $\delta-\gamma n>0$ by hypothesis, if $\gamma>0$, we can choose $\varepsilon>0$ small enough such that the previous inequality fails. This forces $\gamma\le 0$, and the proof is complete. 
\end{proof}

The following result is a generalization of Proposition 7.4. in \cite{CPER}.

\begin{proposition}\label{sharp1}
	Let $1\le p<n$ and $w\in A_r$ with $1\le r\le p$. Suppose that inequality
	\begin{equation*}
		\left( \frac{1}{w(B)} \int_{B} |f(x)-f_B|^{q} w(x)dx\right)^\frac{1}{q}  \le C [w]_{A_r}^\beta  r(B)  \left( \frac{1}{w(B)} \int_{B} | \nabla f(x)|^p w(x)dx\right)^\frac{1}{p}
	\end{equation*}
	holds with some power $\beta$ on the $A_r$ constant. Then $\beta\ge \frac{1}{p}$. 
\end{proposition}

As a consequence of this result, we obtain the sharpness of the exponent $ \frac{1 }{p}$ on the $A_r$ constant in \eqref{eq6}

\begin{proof}
	Let $1<p<n$, let $p\le q <\infty$, and consider the weight $w(x)=|x|^{\delta-n}$ with $0<\delta<nr$, then $w\in A_r$.  Observe that  $w(B)= \frac{\omega_{n-1}}{\delta}$ and  $[w]_{A_r}\simeq _n \frac{n^r}{\delta}\left( \frac{r-1}{nr-\delta}\right)^{r-1}$.  We follow the same ideas as in the previous proof until \eqref{optimality eq1}, where $C_w\sim \left(\frac{n^r}{\delta}\left( \frac{r-1}{nr-\delta}\right)^{r-1}\right)^\beta  $,
	\begin{align*}
		1 \lesssim_{ n} & \frac{n^{r\beta}}{\delta^\beta }\left( \frac{r-1}{nr-\delta}\right)^{(r-1)\beta } \left( 1+ \frac{1}{\left(  1-(2\varepsilon)^\delta \right)^q } \right) \delta^\frac{1}{p} \varepsilon^{\frac{\delta}{p}-\frac{\delta}{q}-1}  \\
		= & n^{r\beta}\left( \frac{r-1}{nr-\delta}\right)^{(r-1)\beta } \left( 1+ \frac{1}{\left(  1-(2\varepsilon)^\delta \right)^q } \right) \frac{\varepsilon^{\frac{\delta}{p}-\frac{\delta}{q}-1}}{\delta^{\beta-\frac{1}{p}} }  .
	\end{align*}
	For a fixed $\varepsilon$, this forces the condition $\beta\ge \frac{1}{p}$. 

\end{proof}

\section*{Acknowledgement} 
I would like to express my deepest gratitude to my advisor, Carlos Pérez, for his guidance and all of the interesting conversations that have led to the writing of this paper. I would also like to thank the anonymous referee for their helpful comments that improved the clarity of the manuscript.

\section*{Conflicts of Interest}
The author has no conflicts of interest to declare.

\section*{Data Availability Statement}
Data sharing not applicable to this article as no datasets were generated or analyzed during the current study.

\bibliographystyle{amsplain}

\end{document}